\documentclass[12pt]{article}

\newcommand{\argm}{\mathop{\rm argmin}\limits}

\large

\begin{document}

%\begin{frontmatter}
%\title{ 

\begin{center}
{\bf Adaptive Optimal Nonparametric Regression and Density Estimation
  based on Fourier-Legendre expansion } \\

\vspace{4mm}

  E. Ostrovsky correspondent author, Y. Zelikov.\\

% address[Ostrovsky] \\ 
 Department of Mathematics, 52900, Bar - Ilan University, Ramat - Gan, 
Israel; \\ 
E - mail: \ galo@list.ru \\
% address[Zelikov] \\ 
  Department of Mathematics, 52900, Bar - Ilan University, Ramat - Gan, 
 Israel; \\
E - mail: zelikoy@macs.biu.ac.il  \\
\end{center}

%\large

\normalsize

\vspace{3mm}

 \begin{abstract}

 Motivated by finance and technical applications, the objective of this paper 
is to consider adaptive estimation of regression and density distribution 
based on Fourier-Legendre expansion, and  construction of confidence intervals
 - also adaptive. The estimators are asymptotically optimal and adaptive in 
the sense that they can adapt to unknown smoothness. \par

\end{abstract}

 % \begin{keyword} 
 Keywords:\par
  Adaptive estimations, regression,  density, martingale,
  confidence interval, Legendre polynomials.
% \end{keyword}

%\end{frontmatter}
\vspace{2mm}

{\it Mathematics Subject Classification 2000.} 41A10(Primary) , 62G07 (Secondary) 

\newpage

\section*{Acknowledgment}
It is the authors pleasure to convey our gratitude to 
M.Lin and V.Fonf (University Beer-Sheva, Israel)
for creative discussions of the problems under consideration, 
to D.Donoho (USA) for sending us his publications and manuscripts, as well as 
A.Pridor, P.Gil and B.Greenstein (Israel) for the 
possibility to practically implement our  estimations.

\newpage

\section{Introduction.} \par

  Among the latest fashions in nonparametric statistics are the so-called 
adaptive estimations (AE), i.e. estimations that use no apriory information 
about the estimated function. Many publications have recently appeared where 
AE are constructed which are optimal in order at a growing number of current 
observations on a continuum of various functional classes (cf. References for 
a list of works on AE, which does not, however, claim to be exhaustive).\par
 In (Polyak B., at al., 1990), (Polyak B. at al., 1992), (Golubev G. at al., 1992) for instance, AE were constructed for the problem of 
estimating regression (R) which are optimal in order on many subspaces of 
space $ L_2 $, and non-adaptive confidence intervals were elaborated on the 
basis of the obtained estimations for the estimated regression function also 
in norm $ L_2 $, which later were somewhat improved in (Golubev at al., 1992).\par

 In  (Efroimovich S., 1985)  AE were 
constructed for problem $(D)$ of estimating distribution density, which are 
optimal on ellipsoids in $L_2 $. \par
 In numerous publications by D. Donoho et al. (Donoho D at al., 1993(1), 1993(2),
1996,  1999(1), 1999(2) )  and in some others AE are 
constructed (and implemented) which are optimal in order on a number of Besov 
spaces. In those papers as well as in  (Golybev G. at al., 1994), (Nussbaum M., 1985),
(Tony Cai at al., 1999), (Lee G., 2003)  diverse orthonormalized 
systems of functions are used to construct AE, 
such as wavelets, wedgelets, unconditional bases, splines, Demmler - Reinsch 
bases, Ridgelets  (Candes E.J., 2003), (Dette H., 2003) etc. \par

 The recent results about kernel estimations in the considered problems see,
for example, (AAD W Van Der Vaart at al., 2003), Allal J., at al.,
2003), (Corinne Berzin at al., 2003).\par

  In  (Ostrovsky E.I., 1996, 1997(1); 1997(2), 
1999) AE were constructed  on the basis of the trigonometric approximation 
theory. \par 

{\bf In this work we construct AE based on the orthogonal polynomial 
expansion series  - the  Legendre polynomials.} \par 
 
{\bf The AE proposed herein feature a speed of convergence  which is optimal 
in order on any regular subspace compactly embedded in space $L_2 $, the 
estimations are universal and very simple in form, which significantly 
facilitates their implementation; finally, we construct exponential adaptive 
confidence intervals (ACI), i.e. such that the tail of the confidence 
probability decreases with exponential speed.}\par

\section{Problem statement. Denotations. Conditions.}

{\bf R.} Regression problem. Let $ f(x), \ x \in [-1,1] $ be an unknown 
function, Riemann-integrable with a square and measured at points of the net 
$ x_i = x_{i,n}= -1+2i/n, \ i = 1,2,\ldots,n; \ n \ge 16 $ with random 
independent centered: $ {\bf E} \xi_i = 0 $ identically distributed 
errors $ \{\xi_i\}: \ y_i = f(x_i) + \xi_i. $ 
It is required to estimate the function $ f(x) $ with the best possible 
precision from the values $ \{y_i\} $.\par
{\bf D.\ Estimation of distribution density.} On the basis of a set of 
independent identically distributed values $\{\xi_i\}, \ \xi_i \in [-1,1], \ 
i = 1,2,\ldots,n $ it is required to estimate their common density $f(x) $ 
(assumed to exist).\par

 It is supposed that all the estimated functions $ f(\cdot) \in L_2[-1,1], $ 
therefore they are expanded in the norm of this space into a Fourier-Legendre 
series in the complete orthonormal system $\{L_j(\cdot)\} $ on the set [-1, 1]: 
$$
f(x) = \sum_{j=0}^{\infty} c_j L_j(x); \ \ c_j = \int_{-1}^1 L_j(x) f(x) dx,
$$
 where $ L_j(\cdot), \ j = 0,1,2,\ldots $ are normalized Legendre's 
polynomials. The Legendre polynomials are given by the Rodrigues formula:
$$ 
P_m(x)={1 \over {2^m m!} } {d^m \over {d x^m}} \bigl [(x^2-1)^m \bigr]
$$
with orthogonal property: 
$$ 
I(k,m) \stackrel{def}{=}
\int_{-1}^1 P_m(x)P_k(x)dx= 2/(2m+1), \ m = k,
$$
 otherwise $ I(k,m) = 0. $  We can define 
$$
L_k(x) = P_k(x) \sqrt{k + 0.5}.
$$
 Let us set $ \rho(N) = \rho(f,N) = \sum_{j=N+1}^{\infty} c_j^2$. 
Evidently $ \lim_{N \to \infty} \rho(N) = 0.$ 
Let us also assume that only the non-trivial {\it infinite-dimensional case}
 will be considered, when an infinite multitude of Fourier coefficients 
$ f $ differs from zero, i.e. $ \forall N \ge 1 \ \Rightarrow \rho(N) > 0.$
 Otherwise our estimations will converge in the sense $ L_2(\Omega \times [-1,1]) $ with speed $1/\sqrt{n}. $ \par
  Moreover, we assume further that (essentially infinite-dimensional case) 

$$
\lim_{N \to \infty} |\rho(N)|/\log N = +\infty. \eqno(2.0).
$$
 In other words, the condition (2.0) means that there exists the constant 
$ q \in (0,1) $ such that for all sufficiently great values $ N  \ \Rightarrow 
 \ \rho(N) \ge q^N. $ \par

 The value $ \rho(N) = \rho(f,N) $ is known and is well studied in the 
{\it approximation theory. } Namely, $ \rho(f,N) = E^2_{N}(f)_2, $ where 
$ E_{N}(f)_p $ is the error of the best approximation of $ f $ 
by the algebraic polynomials of power not exceeding $ N $ in the $ L_p $ 
distance:  for $ g: [-1,1] \to R^1 $ we will  denote 

$$ 
||g||_p = \left(\int_{-1}^1 |g(x)|dx \right)^{1/p}, p \in [1, \infty);  \ 
||g||_{\infty} = \sup_{x \in [-1,1]} |g(x)|, 
$$
 and closely connected with module of continuity of the form

$$ 
\omega_{p,2}(f^{(k)},\delta) = \sup_{h: |h|\le \delta}|f^{(k)}(x+h) - 2 
f^{(k)}(x) + f^{(k)}(x-h)|_p,
$$
for instance:

$$
E_n(f)_p \le C(p,r) \ n^{-r} \  \omega_{p,2} (f^{(r)},1/n).
$$

[DeVore, Lorentz, p. 219-223];  arithmetical operations on the 
arguments of function 
$ f $ and their derivatives are understood as follows: at $ h > 0 \
 x+h = \min(x+h,1),  \ x - h = \max(x-h,-1), \ h > 0. $ \\

Everywhere below condition $ (\gamma1) $ will be considered fulfilled:

$$
(\gamma1): \overline{\lim}_{N \to \infty} \rho(2N)/\rho(N) \stackrel{def}{=} 
\gamma < 1, \eqno(2.1)
$$
	
sometimes stronger conditions $ (\gamma) $ as well:
$$
(\gamma): \exists \lim_{N \to \infty} \rho(2N)/\rho(N) \stackrel{def}{=} 
\gamma < 1;  \eqno(2.2)
$$

$$ 
(\gamma0): \gamma = 0. \eqno(2.3)
$$ 

It is easy to show that it follows from condition (2.0)

$$
\rho(N) \le C N^{-2\beta}, \ 2 \beta \stackrel{def}{=} \log_2(1/\gamma) > 0.
\eqno(2.4)
$$ 

 In the problem (R) it will be assumed that $ \beta > 1/2. $ There are some 
grounds 
to suppose that at $ \beta < 1/2 $ asymptotically optimal AE do not exists in 
the regression problem; for a similarly stated problem this was proved by 
Lepsky (Lepsky O., 1990). \par
  Here and below the symbols 
$ C,C_r $ will denote positive finite constructive constants inessential 
in this context, $ \asymp $ is the usually  symbol, in detail:
$$
A(n) \asymp B(n) \ \Leftrightarrow C_1 \le \liminf_{n \to \infty} A(n)/B(n) 
\le 
$$
$$
\limsup_{n \to \infty} A(n)/B(n) \le C_2, \ \exists C_1, C_2 \in (0,\infty).
$$
 the symbol $ A \sim B $ means that in the given 
concrete passage to the limit $ \lim A/B = 1.$ \par

{\bf Examples.} Denote by $ W(C,\alpha,\beta) $ a class of functions 
$ \{f \} $ such that

$$
\rho(f,N) \sim C N^{-2\beta} (\log N)^{\alpha}, \ \exists C,\beta >0; 
\alpha = const;
$$
 $ W(\alpha,\beta) = \cup_{C>0} W(C,\alpha,\beta);$
$$ 
W(\beta) = W(0,\beta); \ \ W = \cup_{\beta > 0} W(\beta).
$$
 For the class of functions $ W $ the condition $ (\gamma) $ is fulfilled. \par

 Also let us denote
$$
Z(\alpha,\beta) = \{f: \rho(f,N) \sim \alpha \ \beta^N \}, \alpha > 0, \ \beta 
\in (0,1);
$$
and $ Z = \cup_{\alpha > 0; \beta \in (0,1)} Z(\alpha,\beta) $. For the 
functions of class $ Z $ the 
condition $(\gamma0) $ is also fulfilled. Besides, functions of class $ Z $ 
are analytical.\par
 Denote for the problems $ {\bf R, D} $ respectively at $ j < n $ \ \ 
$ \hat{c}_j =  $

$$ 
  (1/n) \sum_{i=1}^n y_j L_j(x_i); \ \ \hat{c}_j = (1/n)
\sum_{i=1}^n L_j(\xi_i);  \eqno(2.5)
$$ 
$ j = 0,1,2,3,\ldots,n-1; $ and for the regression problem
$$
 \ \ B(n,N)= \sum_{k=N+1}^{2N} c_k(n)^2 + \sigma^2 N/n;
$$

$ \sigma^2 = {\bf Var} [\xi_i]; $ for the problem (D) we define 
$ \sigma^2 = 1 $ and
$$
B(n,N) = \sum_{k=N+1}^{2N} c_k^2 +  N/n;
$$
and again for both the considered problems set $ B(n) = $
$$
  \min_{N = 1,2,\ldots, [n/3]} B(n,N), \ N^0 = N^0(n) = 
\argm_{N = 1,2,\ldots,[n/3]} B(n,N);
$$
$$
A(n,N) = \rho(N) + \sigma^2 N/n, \ A(n) = \min_{N = 1,2,\ldots,[n/3]} A(n,N).
$$

 For instance, suppose that $ f \in W(C,\alpha,\beta) $, then $ A(n) \asymp 
n^{-2\beta/(2\beta+1)} (\log n)^{\alpha/(2\beta+1)} $, and in case 
$  f \in Z(\alpha,\beta) \ \Rightarrow A(n) \asymp \log n/n.$ \par
 Our notation should not be surprising, as it follows from condition 
$ (\gamma1) $ that all the introduced functionals $ \{B(n,N)\}, \ \{B(n)\} $ arising from different problems are mutually  $ \asymp $  equivalent. Besides, for the same reasons	
$$
A(n,N) \asymp B(n,N); \ A(n) \asymp B(n).
$$

 Apart from that it is clear that in the regression problem conditions must 
be imposed not only on the estimated function, but on the measurement errors 
$ \xi_i $ too.  We will consider here only the so-called 
{\it exponential} level. Indeed, we assume that in the regression problem the 
following condition is satisfied: 

$$
(Rq): \exists q,Q \in (0,\infty), \Rightarrow {\bf P} (|\xi_i| > x) \le 
\exp \left(-(x/Q)^q \right), x > 0.
$$
 
 The so-called classical projective estimates was introduced by 
N.N.Tchentsov [Tchentsov N.N., 1972, p. 286] (for the trigonometrical 
system instead considered here Legendre's polynomials $ L_k(\cdot) ) $
will be considered as an estimates of the function $ f $:
$$
f(n,N,x) = \sum_{j=0}^N \hat{c}_j L_j(x). \eqno(2.6)
$$
 Since, as is shown by Tchentsov, 
$ {\bf E}||f(n,N,\cdot) - f(\cdot)||^2 \asymp B(n,N), $ the selection of the 
number of summands $ N $ optimal by order in the sense of 
$ L_2(\Omega) \times L_2[-1,1] $ is given by the expression $ N = N^0(n) $
with the speed of convergence  $ f(n,N^0,\cdot) \to f(\cdot) $
in the above-mentioned sense is $ \sqrt{A(n)}. $ I. A. Ibragimov and R. Z. 
Khasminsky (Ibragimov I., Khasminsky R., 1982)
 proved that no faster convergence exists on the regular classes of 
functions given by the value $ \sqrt{A(n)}. $ \par
 However, the value $ \rho(f,N) $ or at least its order as $ N \to \infty $ are 
practically unknown as a rule. Below the adaptive estimation of $ f $ will be 
studied based only on the observations $ \{\xi_i\} $ and using no apriory 
information regarding $ f $, and yet possessing the optimal speed of 
convergence at apparently weak restrictions. Set 

$$ 	
\tau(N) = \tau(n,N) \stackrel{def}{=} \sum_{k=N+1}^{2N} \hat{c}^2_k, \
 \ N(n) \stackrel{def}{=} \argm_{N \in (1, [n/3])} \tau(n,N),  \eqno(2.7)
$$ 

$$
\tau^*(n) = \min_{N \in (1,[n/3])} \tau(n,N), \ 
$$
{\it Our adaptive estimations $ \hat{f} $ in both considering problems 
 have a universal view:}

$$ 
\hat{f} = f(n,N(n),x) = \sum_{0=1}^{N(n)} \hat{c}_j L_j(x).  \eqno(2.8)
$$ 

 In case of a non-unique number of summands $ N(n) $ in (2.7) we choose the 
largest. Below the value $ N $ will always be {\it arbitrary } non-random 
integer number in the set of integers numbers 
of the segment $ 1,2,\ldots,[n/3] $ and $ N(n) $ is the {\it random }
variable defined in (2.7). \par

 Note, that by using the  Fast Legendre Transform technique described by D. Potts et al. (D. Potts et al., 1998), the amount of elementary operations for  $ \hat{f} $ calculation 
is  $O( n \log  n)$, likewise in Fast Fourier Transform and in Fast Wavelet Transform. \par

 Before proceeding to formulations and proofs let us clarify informally our 
idea for choosing $ N(n). $ It is easy to find by direct calculation 
for the regression problem (and analogously for the problem D) 
that the coefficients estimations $ \hat{c_k} $ have a view:

$$
\hat{c_k} = c_k(n) + n^{-1/2} \theta_k(n),
$$

where $ as n \to \infty $

$$
c_k(n) = n^{-1} \sum_{i=1}^n f(x_i) L_k(x_i) \to \int_{-1}^1 f(x) L_k(x) dx =
 c_k;
$$

$$
\theta_k(n) = n^{-1/2} \sum_{i=1}^n \xi_i L_k(x_i).
$$

 It follows from the multidimensional CLT that the variables  
$ \ \{ \theta_k(n) \} $ as $ \ n \to \infty $ are asymptotically Gaussian 
distributed and independent:

$$
{\bf Var} [\theta_k(n)] = n^{-1} \sum_{i=1}^n \sigma^2 L^2_k(x_i) 
\to \sigma^2 \ \int_{-1}^1 L^2_k(x) dx = \sigma^2;
$$

$$
{\bf E} \theta_k(n) \theta_l(n)  = \sigma^2 n^{-1} \sum_{i=1}^n L_k(x_i) 
L_l(x_i) \to \sigma^2 \ \int_{-1}^1 L_k)x) L_l(x) \ dx = 0, \ k \ne l.
$$
 Therefore, the variables $  \{ \theta_k(n) \} $
are asymptotically independent and have approximately the normal distribution: 
$$
Law(\hat{c_k} ) \asymp N(c_k,\sigma^2/n),
$$
or equally 
$$ 
\hat{c_k} = c_k + \sigma \epsilon_k/\sqrt{n}, \ Law(\epsilon_k) \asymp N(0,1)
$$
and also $ \{\epsilon_k\} $ are asymptotically independent. Therefore,
$ \tau(n,N) \asymp $

$$
\sum_{k=N + 1}^{2N} c_k^2 + 2 \ n^{-1/2} \ \sigma 
\sum_{k = N + 1}^{2N} c_k \epsilon_k + \sigma^2 \ n^{-1} 
\sum_{k = N + 1}^{2N} \epsilon_k^2;
$$

$$ 	
{\bf E} \tau(n,N) \asymp B(n,N), \ \ {\bf Var} [\tau(n,N) ]
\asymp B(n,N)/n,  \eqno(2.9)
$$ 
and therefore 

$$ 
N \to \infty, N/n \to 0 \ \Rightarrow \sqrt{ {\bf Var} [\tau(n,N)]}/
{\bf E}\tau(n,N) \to 0.
$$
 Note that in the case of the regression problem the condition 
$$
\beta > 1/2  \eqno(2.10)
$$
is essential which is common in statistical research (Polyak B. at al., 1990, 
1992),(Lepsky O., 1990).  We will assume in the problem (R) that the 
condition  (2.10) is satisfied.\par
 It follows from (2.9) that there are some grounds to conclude
$$
\tau(n,N) \stackrel{a.s}{\asymp} {\bf E} \tau(n,N) \asymp A(n,N)
$$
and therefore
$$
N(n) = \argm_{N \le n/3} \tau(n,N) \sim \argm_{N \le n/3} {\bf E} \tau(n,N) 
= N^0(n).
$$
 Also note that the number of summands $ N(n) $ proposed by us is a random 
variable (!) and that estimation (2.8) is non-linear by the totality of 
empirical Fourier coefficients $ \{\hat{c}_j\}. $ \par

\vspace{3mm}

\section{ Formulation of the main results.}\par
 Further we will investigate the exactness of our adaptive estimations 
in the $ L_2 $ sense in our interval [-1,1] and will write as usually

$$
||f-g||^2 = \int_{-1}^1 (f(x) - g(x))^2 \ dx.
$$
 Define also for the problem (R) 
$$
r = r(q) = 2q/(q + 4), q \in (0,2); \ q \ge 2 \ \Rightarrow r = q/(q+1),
 \eqno(3.0)
$$
 and $ r = 1 $ for the problem (D). \par

{\bf Theorem R(q).} {\it Under the conditions} $(Rq), (\gamma1),$
 {\it in the problem } {\bf R } {\it we propose that there exists 
  a constant } 
   $ K_R = K_R(q,\gamma) \in (0, \infty) $ {\it such that for the 
 variable }
 
 $$
 \zeta_R = \zeta_R(n) = B^{-1}(n) \ Q^{-2} \ ||\hat{f} - f||^2  - Q^{-2} \ K_R
 $$
{\it  the following inequality holds:}
 
 $$
 {\bf P}(|\zeta_R| > u) \le 2 \exp \left(- C \ u^{r/2} \ 
 (n \ A(n))^{r/2}/\log\log n \right), \ u > 1. \eqno(3.1)
 $$

  (See in comparison (Bobrov P. at al., 1997); here
 the exponent indices are significantly decreased.)\par

 {\bf Theorem D(q).} {\it Under the condition} $ (\gamma1),$
 {\it in the problem } {\bf D } {\it we propose that there exists a constant}
 $ K_D = K_D(q,\gamma) $ {\it such that for the variable }

 $$
 \zeta_D = \zeta_D(n) = B^{-1}(n) ||\hat{f} - f||^2  - K_D
 $$
 {\it the following inequality holds:}
 
 $$
 {\bf P}(|\zeta_D| > u) \le 2 \exp \left(- u^{1/2} \ (n \ A(n))^{1/2}
  /\log\log n \right), \ u > 1. \eqno(3.2)
 $$
 
  This result improves the one for the Fourier approximation of
(Ostrovsky E., Sirota L., 2004). \par

{\bf Theorem $ (Rq) \ a.s.$ } {\it If in the problem (R) under condition 
$ (Rq) $ for arbitrary $ \varepsilon > 0 $ the series}

$$ 
\sum_{ \{n > 16 \} } {\bf P}_n(\varepsilon) \stackrel{def}{= }
\sum_{ \{n > 16 \} } \exp \left(- \varepsilon \frac{(n A(n))^{r/2} }
{\log \log n} \right) < \infty, \eqno(3.3)
$$ 

{\it converges, then in the sense of convergence  with probability one} 

$$
 \lim_{n \to \infty} \tau^*(n)/B(n) = 1 \eqno(3.4 a)
$$
{\it and }

$$
\overline{\lim}_{n \to \infty} || \hat{f} - f||^2/B(n) \le K_R. \eqno(3.4.b)
$$

 Let us make another additional assumption (v) with regard to the class of 
estimated functions $ \{f\}. $ Denote 

$$	
H_n(v) = \inf_{N \in [N^0/v, N^0 v]} \frac{B(n,N)}{B(n)}, \ v = const > 1;
$$

$$ 
 {\bf P}^{(N)}(n,v) = \exp \left( - (\sqrt(H_n(v)) - 1)^{r/2} (n B(n))^{r/2}
/ \log \log n \right). 
$$ 
(At $ v \ge N^0 $ the left interval is absent, at $ v \ge n/(3N^0) $ 
thr right interval is absent.) \par
 {\it Condition } $(v): $

$$
\forall v > 1 \ \Rightarrow \sum_{ n \ge 16} {\bf P}^{(N)}(n,v) < \infty.
$$

 The classes of functions satisfying conditions $ (\gamma1) $ and $(v) $ will 
be called regular. Classes $ W $ and $ Z $ are regular.\par

{\bf Corollary 1.} {\it If in addition to the conditions of theorem  R(q) the 
condition $ (v) $ holds, then a.e.}
$$
\lim_{n \to \infty} N(n)/N^0(n) = 1. \eqno(3.5) 
$$

{\bf Theorem (D) a.s.} {\it Let for problem $ (D) $ besides the 
above-formulated assumptions, condition (3.3) also be fulfilled with 
$ r/2 $ replaced by 1/2. Then the factual convergences of (3.4 a) and 
(3.4.b) are asserted here as well.}\par
{\bf  Corollary 2.}
{\it Analogously if in addition to the our the condition $ (v) $ holds, 
then also with probability one }
$$
\lim_{n \to \infty} N(n)/N^0(n) = 1. 
$$

\vspace{3mm}

{\bf 4. Proofs. }\par
\vspace{2mm}
  The proofs of the theorem (Rq) and (Dq) are similar to 
 proofs of the our result in (Ostrovsky E., Sirota L., 2004) for trigonometrically 
 approximation for $ f(x); $ we will use the known properties of Legendre's 
 polynomials (Kallaev, 1970), (Szeg\"o, 1959). For instance, 

$$
 \sup_k \sup_{x \in [-1,1]} |L_k(x)| \left(1  -  x^2 \right)^{1/4} 
< \infty.
$$
  Instead of the semi invariant estimations for polynomials from independent 
 random variables (Saulis, Statuliavitchius, 1989) we will use the modern  
 estimations for polynomial martingales   (Hall, Heyde, 1980) [pp. 115 - 120],(Ostrovsky. E, 2004). \par 
  First of all we consider the problem of regression (R).\par
  We will assume without loss of generality $ Q = 1. $ \par

 {\sc STEP 1.} Let us write the exact expression for the important variables.
 Introduce the notation:\\

 $$
  \delta_k(n) \stackrel{def}{=} c_k(n) - c_k.
 $$

 We can write:  
  
   $$ 
  \hat{c}_k(n) = c_k + \delta_k(n) + n^{-1} \sum_{i=1}^n \xi_i L_k(x_i),
   $$

  $$
  (\hat{c}_k(n))^2 =  c^2_k + \delta^2_k(n) + \sigma^2 n^{-2} \sum_{i=1}^n 
   L^2_k(x_i) +
  $$
  
  $$
  2 c_k \delta_k(n) + 2 n^{-1} \sum_{i=1}^n  c_k \xi_i L_k(x_i) + 2 n^{-1} 
  \sum_{i=1}^n \delta_k(n) \xi_i L_k(x_i) + 
  $$
  
   $$
    n^{-2} \sum_{i=1}^n \left( \xi^2_i - \sigma^2 \right) L^2_k(x_i) +
    n^{-2} \sum \sum_{i \ne j} \xi_i \xi_j L_k(x_i) L_k(x_j).
   $$

  We have for the variables $ \tau(N,n) $ (and further for the variables
 $ \Delta^2 = \Delta^2(N,n) = ||\hat{f} - f||^2 \ ): $  $ \tau(n,N) = $

$$
\left[ \sum_{k=N+1}^{2N} c^2_k + 2 \sum_{k=N+1}^{2N} c_k \delta_k(n) +
\sum_{k=N+1}^{2N} \delta^2_k(n) + \sigma^2 n^{-1} \sum_{k=N+1}^{2N} 
n^{-1} \sum_{i=1}^n L^2_k(x_i) \right] + 
$$
 
 $$
 \left[ 2 n^{-1} \sum_{i=1}^n \xi_i \sum_{k= N+1}^{2N} \delta_k(n) L_k(x_i)
 + 2 n^{-1} \sum_{i=1}^n \xi_i \sum_{k=N+1}^{2N} c_k L_k(x_i) \right] + \tau_2,
 $$

 $$
\tau_2 = \left[ n^{-1} \sum_{i=1}^n \left(\xi^2_i - \sigma^2 \right) n^{-1} 
\sum_{k=N+1}^{2N} L^2_k(x_i) \right] +
 $$

 $$
 \left[ 2 n^{-1} \sum \sum_{1 \le i < j \le n} \xi_i \xi_j \ n^{-1} 
 \sum_{k = N+1}^{2N} L_k(x_i) L_k(x_j) \right], \eqno(4.1)
 $$

where $ \tau = \tau_0 + \tau_1 + \tau_2; \ \tau_m = \tau_m(n,N), \  \tau_0 $  is the deterministic  part of $ \tau:$
$ {\bf E} \tau = \tau_0 \sim B(n,N), \ \tau_1 $  is the linear combination of  
$ \{ \xi_i \}, \ \tau_2 $ is the bilinear combination of $ \{ \xi_i \}$.\par
 
 It is easy to verify using the known properties of Legendre's 
polynomials that $ \tau_1 \asymp B(n,N) $ and that

$$
{\bf Var} [\tau_1] \le C B(n,N)/n, \ {\bf Var}[\tau_2] \le  C B(n,N)/n.
$$

{\sc STEP 2.} Note that the sequences $ \eta_1(n) = \sum_{i=1}^n b(i) \xi(i),
 \ \eta_2(n) = \sum_{i=1}^n  b(i)(\xi^2_i - \sigma^2) $ and 

$$
\eta_3(n) = \sum \sum_{1 \le i < j \le n} b(i,j) \xi_i \xi_j, 
$$

where $ \{ b(i) \}, \{b(i,j) \} $ are a non-random sequences, with the second 
component $ F(n) = \sigma \left( \{ \xi_i \}, i = 1,2,\ldots,n \right), $ i.e.
$ \{ \eta_s(n), F(n) \}, s = 1,2,3;  \ \{ F(n) \} $ is
the natural sequence (flow) of sigma-algebras,  are martingales. \par
 It follows from the main result of paper (Ostrovsky E., Sirota L., 2004), 
devoted to the exponential and moment estimations for martingale distributions, that

$$
\sup_{n \ge 16} {\bf P}( | \tau_{1,2} | / \sqrt{ {\bf Var} \left[\tau_{1,2}
\right]} > x) \le \exp \left( - C x^r \right), \ x > 0. \eqno(4.2)
$$

{\sc STEP 3.}
  Now we intend to use on the basis of inequality (4.2) the Law of Iterated 
 Logarithm (LIL) for the martingales (Hall, Heyde, 1980) [pp. 115 - 121] in the 
more convenient for us form (Ostrovsky, 1999) [pp. 79 - 83]. Namely, if we denote

$$
\nu(n,N) = ( \tau(n,N) - {\bf E} \tau(n,N) )/\sqrt{B(n,N)/n},
$$
then $ {\bf E} \nu(n,N) = 0 $ and if we denote

$$
\zeta \stackrel{def}{=} \sup_{n \ge 16} \sup_{ N \in [1, n/3]} |\nu(n,N)| /
\log \log n,
$$

then $ \zeta < \infty $ a.e. and for the random variable $ \zeta $
we have  for all positive values $ x, x > 0 $ the tail inequality

$$
{\bf P} ( \zeta > x) \le \exp \left( - C(q) x^r \right). \eqno(4.3)
$$
 The inequality (4.3) may be rewritten as follows:

$$
\tau(n,N) = {\bf E} \tau(n,N)  + \zeta(n,N) \ \log \log n \ \sqrt{B(n,N)/n},
\eqno(4.4)
$$

where $ \zeta =  \sup_{n,N} |\zeta(n,N)| $ satisfies the inequality (4.3), 
$ {\bf E} \tau(n,N) = B(n,N)(1 + \theta(n)), \ \theta(n) $ is non-random and
$ \lim_{n \to \infty} \theta(n) = 0. $  \par

\vspace{2mm}

{\sc STEP 4.} Let $ M $ be some subset of an integer segment 
$ S = [1,2,\ldots,n], $ $ \overline{M} = S \setminus M, \ 
\pi(M) \stackrel{def}{=} {\bf P}(N(n) \in M), $ and assume that
$$
v =v(n,M) \stackrel{def}{=} \inf_{N \in M} B(n,N)/B(n) > 1.
$$
Then under conditions $ (\gamma) $ and $ (Rq) $

$$ 
\pi(M) \le 2\exp \left(- C \  \left[ ( \sqrt{v} -1) \ n B(n) \right]^{r/2} /
 (\log \log n) \right). \eqno(4.5)
$$ 

{\bf Proof.} We obtain for the case of $ (Rq) $, denoting $ \overline{\nu} = 
\max_{N \in S} |\nu(n,N)|: $  
$$
\pi(M) = {\bf P}(N(n) \in M) = {\bf P} (\min_{N \in \overline{M}} \tau(n,N) >
\min_{N \in M} \tau(n,N)) =
$$

$$
 {\bf P} \left( \min_{N \in \overline{M}} ( B(n,N) + \sqrt{B(n,N)/n}
 \ (\log \log n)^{1/r} \ \nu(n,N)) \right) >
$$

$$
 \min_{N \in M} \left( B(n,N) + \sqrt{B(n,N)/n} \ (\log \log n)^{1/r} 
\ \nu(n,N) \right)  \le 
$$

$$
 {\bf P} ( B(n) + \sqrt{B(n)/n} \ (\log \log n)^{1/r} \ \overline{\nu} >
 v B(n) - \sqrt{v B(n)/n} \ (\log \log n )^{1/r} \ \overline{\nu} ).
$$

 We find solving the inequality under the probability symbol relative to 

$$
 \overline{\nu}: \ \pi(M) \le  {\bf P} \left( \overline{\nu} \ (1 + \sqrt{v}) \sqrt{B(n)/n} \ (\log \log n)^{1/r} \ge (v-1) B(n) \right) \le 
$$

$$
 {\bf P}\left( \overline{\nu} \ge \frac{v-1}{ \sqrt{v}+ 1}  \ 
\frac{ \sqrt{n B(n)} } {(\log \log n)^{1/r}} \right) =
 {\bf P} \left( \overline{\nu} \ge   (\sqrt{v} - 1) \frac{\sqrt{n B(n)}}
{(\log \log n)^{1/r} } \right). \eqno(4.6)
$$
 Using our estimations (4.6) for $ \tau, $ we arrive to the assertion (4.5).\par
 Note that under our condition (2.0) 

$$
n B(n) > C \ \log n, \ C_1 \in (0,\infty),
$$
therefore under our conditions  for all values $ v, \ v > 1 $ 
$$
 \lim_{n \to \infty} \pi(M) = 0.
$$
 If in addition for any $ \varepsilon > 0 $ the series 
$ \sum_n {\bf P}_n(\varepsilon) $ converges, the assertions (3.4.a),(3.4.b) 
and  corollaries 1,2
to be proved follows from the  lemma of Borel-Cantelli. \par

 The rest is proved analogously if it is taken into account that 
$ N \ge N^0(n) (1 + \varepsilon), \ \varepsilon \in (0,1] $ and condition 
$ (v) $ lead to the inequality $ B(n,N) \ge (1 + C \varepsilon^2) B(n), $ 
$ \varepsilon \in (0,1); $ this completes the proof.\par
 Analogously we can prove the theorem {\bf (Rq)a.s}, on the basis of 
inequality: 
$$
{\bf P}_n (\varepsilon) \le 
\exp \left( - C \varepsilon^{r/2} \ (n A(n)^{r/2})/ (\log \log n) \right).
$$

 Note in addition that  at $ v > 2 $

$$
\tau^*(n) \le B(n)(1 + \theta(n)) + \overline{\nu} \ \log \log n \ 
\sqrt{B(n)/n},
$$

therefore we have for sufficiently great values $ n: $

$$
{\bf P} ( |\tau^*(n)/(B(n)(1 + \theta(n)))| > v ) \le
\exp \left( - C \varepsilon^{r/2} \ (n A(n)^{r/2})/ (\log \log n) \right)
\eqno(4.7)
$$

and analogously

$$
{\bf P} ( |\tau^*(n)/(B(n)(1 + \theta(n)))| < 1/v ) \le
\exp \left( - C \varepsilon^{r/2} \ (n A(n)^{r/2})/ (\log \log n) \right).
\eqno(4.8)
$$

\vspace{3mm}

{\sc STEP 5.} Let us consider here the main variables $ \Delta^2. $  We have:
$ \Delta^2 = $
$$
= \left[ \sum_{k=N(n)+1}^{\infty} c^2_k + \sigma^2 n^{-1} \sum_{k=0}^{N(n)} 
n^{-1} \sum_{i=1}^n L^2_k(x_i) +
 \sum_{k=0}^{N(n)} \delta^2_k(n) + 2 \sum_{k=0}^{N(n)} c_k \delta_k(n) \right] +
 $$

 $$
\left[ 2 n^{-1} \sum_{i=1}^n \xi_i \sum_{k= 0}^{N} \delta_k(n) L_k(x_i)
 + 2 n^{-1} \sum_{i=1}^n \xi_i \sum_{k=0}^{N} c_k L_k(x_i) \right] + 
 $$

 $$
 \left[ n^{-1} \sum_{i=1}^n \left(\xi^2_i - \sigma^2 \right) n^{-1} 
\sum_{k=0}^{N(n)} L^2_k(x_i) +
2 n^{-1} \sum \sum_{1 \le i < j \le n} \xi_i \xi_j \ n^{-1} \sum_{k = 0}^{N(n)} 
 L_k(x_i) L_k(x_j) \right] =
 $$

 $ \Delta_0 + \Delta_1 + \Delta_2. $ We have analogously to the investigation 
of the expression  for $ \tau(n,N) $ using the condition  $ \gamma $
 for sufficiently large  values  $ n \ge n_0 > 2: $
$$
||\hat{f} - f||^2/B(n) \le C A(n,N(n))/B(n) + \Psi_3(N(n)) /B(n) \le 
$$
$$
C(1-\gamma)^{-1} \tau^*(n)/B(n) + \Psi_3(N(n))/B(n) =
$$
$$
 C (1-\gamma)^{-1}  + (\tau^*(n) / B(n) - 1) + \Psi_3(N(n)) /B(n),
$$
where, as can easily be seen, $ \Psi_3(N) = \Delta_1 + \Delta_2.  $ \par

 Then we will use the elementary inequality 
$ {\bf P}({\bf A}) \le {\bf P}( {\bf AB}) + {\bf P}({\bf \overline{B}}),$
 in which $ {\bf A,B} $ are events. Setting $ {\bf A} = $
$$
 \{ ||\hat{f} - f||^2/B(n) - C/(1 - \gamma) > u\}, \ 
{\bf B} = \{ 1/v \le \tau^*(n)/ B(n) \le v \},
$$

 we have  at $ v \in (2, u - C): $ 
$$
{\bf P}_0 \stackrel{def}{=}{\bf P} ({\bf AB} ) \le 
{\bf P}(v+ \max_{ N: 1/v < \tau^*/B(n) < v } |\Psi_3(N)|/B(n) > u).
$$
We find analogously to the (Ostrovsky E., Sirota L., 2004): 

$$
{\bf P}(v+ \max_{ N^0/v < N < v N^0} |\Psi_3(N)|/B(n) > u) \le
$$

$$
\exp \left(-C \frac{(u-v)^r}{v^{r/2}} \ \frac{((n A(n))^{r/2}}
{\log \log n} \right).
$$
 Thus, $ {\bf P}(A) \le $

$$
\exp \left(-C \frac{(u-v)^r}{v^{r/2}} \ \frac{((n A(n))^{r/2}}
{\log \log n} \right) + \exp \left(-C v^r \ \frac{((n A(n))^{r/2}}
{\log \log n} \right).
$$

Choosing $ v = 0.5u $ for sufficiently great values $ u, \ u \ge C, \ $ we arrive to the assertion of theorem R(q).\par

 {\bf Remark 1.} Let us note, and use it below, a slight difference in the 
behaviors of the values $ \tau(n,N) $ and $ N(n) $ which consists in the 
peculiarity of condition $ (v). $ At $ v > 1 $ we have (under the same 
conditions $(Rq),(v): $
$$
\max \left( {\bf P} \left(\frac{N(n)}{N^0(n)} \le \frac{1}{v} \right),
{\bf P} \left( \frac{N(n)}{N^0(n)} > v \right) \right) \le \exp \left( - C v^r \frac{(n A(n))^{r/2} }{(\log \log n ) } \right).
$$
 An analogous estimation for the probability $ {\bf P}(\tau^*(n)/B(n) > v) $ 
holds even without condition $ (v).$ \par
{\bf Remark 2.} The {\it consistency } of the proposed estimations in the 
above-mentioned sense under all the introduced conditions, including (v),
it follows from the assertions already proved. Indeed, since 
$$ 
A(n) \le A(n,[\sqrt{n}]) \le C n^{-1/2} + \rho([\sqrt{n}]) \to 0, 
$$ 
 then $ N^0(n) \to \infty, \ N^0(n) / n \to 0, $ because otherwise the value 
$$ 
A(n) = A(n,N^0(n)) \asymp N^0(n) / n+ \rho(N^0)
$$ 
would not tend to zero. \par
 Since $ N(n)/N^0(n) \to 1, $ then $ N(n) \to \infty $ and analogously 
$ N(n)/n \to 0, $ which proves the consistency of $ \hat{f}.$ \par

{\bf We proceed now to the problem of estimating density (D).}  Here 

$$
c_k = \int_{-1}^1 L_k(x) \ f(x) dx = {\bf E} L_k(\xi_i), \
\hat{c}_k = n^{-1} \sum_{i=1}^n L_k(\xi_i), 
$$

$$
 \hat{f}(x) = \sum_{k=0}^{N(n)} \hat{c}_k \ L_k(x) =
n^{-1} \sum_{i=1}^n \sum_{k=0}^{N(n)} L_k(\xi_i) \ L_k(x).
$$

The functional $ \tau $  may be written as $ \tau(n,N) = $

$$
 n^{-1} \sum_{i=1}^n n^{-1} \sum_{k=N+1}^{2N} L^2_k(\xi_i) +
2 n^{-2} \sum_{1 \le i < j \le n} \ \sum_{k=N+1}^{2 N} L_k(\xi_i) L_k(\xi_j).
$$

Let us denote 

$$
G_k(x,y) =  \sum_{m=0}^k L_m(x) L_m(y).
$$
It is known (Bateman H., Erdelyi A., 1953)  [chapter 10, section 10], that if 
$ x \neq y, $ then 

$$
G_k(x,y) = (k+1) \left[P_{k+1}(x) P_k(y) - P_k(x) P_{k+1}(y) \right]]/(x - y)
$$

and 

$$
G_k(y,y) = (k+1) \left[P_k(y) P^/_{k+1}(y) - P_{k+1}(y) P^/_k(y) \right].
$$

Also 

$$
\left(1-x^2  \right) P^/_k(x) = k \left[P_{k-1}(x) - x P_k(x) \right].
$$

 Therefore, we have in the considered problem (D) 

$$
\tau(n,N) = n^{-1} \sum_{i=1}^n n^{-1} \left( G_{2N}(\xi_i) - 
G_N(\xi_i) \right) +
$$

$$
2 n^{-2} \sum \sum_{1 \le i < j \le n} \left(G_{2N}(\xi_i, \xi_j) -  
G_N(\xi_i, \xi_j) \right) = 
$$
$ \tau_0 + \tau_1 + \tau_2, \ {\bf E} \tau = \tau_0 \sim B(n,N), $ and
the second (and the first) expression for the $ \tau, $ i.e. 
$ \tau_2 $ is the so-called $ U $ statistics. \par
 We find by direct calculation  (as in the case of problem R):

$$
{\bf Var} [\tau_1] \le C B(n,N)/n, \ {\bf Var}[\tau_2] \le  C B(n,N)/n.
$$

 Recall that the $U$ – statistic with correspondent sequence of sigma-algebras 
is also a martingale. Using the exponential boundaries for the martingale 
distribution, (Ostrovsky E., 2004), (Korolyuk B.S., Borovskich Yu., 1993) etc.,
we obtain:

$$
\sup_{n \ge 16} {\bf P}( | \tau_{1,2} | / \sqrt{ {\bf Var} \left[\tau_{1,2}
\right]} > x) \le \exp \left( - C x \right), \ x > 0, 
$$

$$
{\bf P} ( \zeta > x) \le \exp \left( - C \ x \right),
$$
where $ \zeta =  \sup_{n,N} |\zeta(n,N)|, $

$$
\tau(n,N) \stackrel{def}{=} {\bf E} \tau(n,N)  + \zeta(n,N) \ \log \log n \ \sqrt{B(n,N)/n}.
$$

 Repeating the considerations of (Ostrovsky E., Sirota L., 2004] we complete the proof. \par

{\bf 5. Adaptive confidence intervals (ACI).} Let us now describe the use of our 
results for the construction of ACI. Note first of all that the probability 
$$
 {\bf P}_f(u)  = {\bf P}(||\hat{f} - f||^2 > u)
$$ 
with rather weak conditions in all the 
considered problems permits estimation of the form

$$ 
{\bf P}_f(u) \le 2 \exp \left(- \varphi(C,n, B(n)) u^{r/2}) \right)
\stackrel{def}{=} {\bf P^+_f}(u),
\ u > C_1.  \eqno(5.1)
$$ 
 As it is proved above, the variables $ B(n), C, C_1 $ have respective 
consistent estimates, for example,
$$
 B(n) \approx \min_{N \le n/3} \tau(n,N) = \tau^*(n).
$$
 The values $ C, C_1 $ also depends on $ \gamma $ and on the constants 
$ C_j $ appearing in the definition of conditions $ (\gamma), (v) $. With 
very weak conditions 
they can also be estimated consistently by the sampling in the following way. 
Set $ M = M(n) = \left[\exp(\sqrt{\log n}) \right] $; then, if conditions 
$ (\gamma), (v) $ are fulfilled, a system of an asymptotic equalities can be 
written:
$$
\tau(M) - \sigma_s M/n \sim (1-\gamma)\rho(M);
$$
$$
\tau(2M) - 2 \sigma_s M/n  \sim \gamma (1-\gamma) \rho(M);
$$
$$
\tau(4M) - 4 \sigma_s M/n \sim \gamma^2 (1-\gamma) \rho(M),
$$
where the symbol $ s $ denotes the number of problem.\par
 Solving this system, we find the consistent $ (mod \ \ {\bf P}) $ estimate 
of $ \gamma: $
$$
\hat{\gamma} = \frac{\tau(4M) - 2 \tau(2M)}{\tau(2M) - 2 \tau(M)}.
$$
(The parameter $ \sigma_s $ can also be estimated consistently, but that 
is not necessary for us). The constants $ C_1, C $ also can be consistent 
determined. \par

 Substituting the obtained estimates of all the parameters 
into (5.1), we get to the estimate of the confidence probability

$$ 
 {\bf P^+_f} (u) \le 2 \exp \left( - \phi(C(\hat{\gamma}, \hat{C}_1,\hat{C} ), 
n,\tau^*(n)) \ u^{r/2} \right)
\stackrel{def}{=} \hat{ {\bf P}}_f(u).\eqno(5.2)
$$ 
then, equating the right-hand part of (5.2) of the unreliability of the 
confidence interval $ \delta $ to, say, the magnitude 0.05 or 0.01, 
we calculate $ u = u(\delta) $ from the relation
$$
\hat{ {\bf P}}_f(u(\delta)) =\delta
$$
and obtain approximately the {\it adaptive confidence interval} for $ f $ 
reliability $ 1 - \delta $ of the form

$$ 
|| \hat{f} - f||^2 \le u(\delta) \min_{N \le n/3} \tau(n,N).\eqno(5.3)
$$
But for a rough estimate of the error from replacing $ f $ by $ \hat{f} $ the 
following quite simple method can be recommended. Since

$$ 
\frac{||\hat{f} - f||^2}{ B(n)} = \frac{A(n,N(n))}{B(n)} + 
\frac{\Psi_3(N(n))}{B(n)},\eqno(5.4)
$$ 

and the second term in the right-hand part of (5.4) a.s. tends to 
zero, while the first term, if conditions $ (\gamma), (v) $ are fulfilled, 
has 
$ 1/(1 - \gamma) $ as its limit, we thus prove the following assertion 
apparently well known to specialists in nonparametric statistics for 
non-adaptive estimation: \newline

{\bf Theorem c.i.} {\it If the following conditions are fulfilled in our 
problems: $ $ in the problem $ R \ (Rq), (\gamma), (v) $ or $ (\gamma), (v) $ in 
problems $ D, \ S, $ then}

$$
\overline{\lim}_{n \to \infty} ||\hat{f} - f||^2 /B(n) \le 1/(1 - \gamma).\eqno(5.5)
$$ 

 In order to construct an adaptive confidence interval assertion (5.5) can be 
reformulated as follows. {\it With probability tending to 1 as } 
$ n \to \infty $
$$
|| \hat{f} - f||^2 \le B(n)/(1 - \gamma), \eqno(5.6)
$$
and ACI is constructed by replacing the values $ B(n), \gamma $ by their 
consistent estimates:
$$
|| \hat{f} - f ||^2 \le \tau^*(n) \ \frac{ \tau(2M) - 2 \tau(M)}{ 3 \tau(2M) - 2 
\tau(M)  - \tau(4M) }. \eqno(5.7)
$$
A more exact result will be obtained by taking into account the following term 
of the expansion of the value $ || \hat{f} - f ||^2: $
$$
\frac{ || \hat{f} - f ||^2}{ B(n) } \le \frac{1}{ 1 - \gamma} + \frac{\zeta}
{\sqrt{N^0(n)}} (1 + \epsilon_n),
$$
where $ \epsilon_n \to 0; \ {\bf P}(|\zeta| > u) \le 2 \exp( - C u^{r/2}) $ and 
$ C $ no longer depends on $ n $. Equating the probability 
$ {\bf P}(|\zeta| > u) $, more exactly its estimate 
$ 2 \exp( - C u^{r/2} ) $ to the value $ \delta, \ \delta \approx 0+, $ we will 
easily find $ u = u(\delta) $ and construct an approximate ACI with reliability 
$ \approx 1 - \delta $ of the form
$$
|| \hat{f} - f ||^2 \le \frac{ \tau^*(n)}{1 - \hat{ \gamma}} + \tau^*(n) 
u(\delta).
$$
Closer consideration reveals an effect that somewhat reduces the exactness of 
ACI. Let (as is true in all the three considered problems under the formulated 
assumptions)

$$
{\bf P} \left( ||\hat{f} - f||^2 / B(n) > u \right) \le \exp( 
- \phi(C_1 u)), \ \phi(u) = \phi(n,u),
$$
$$
{\bf P}\left( \tau^*(n)/B(n) < 1/u \right) \le 
\exp( - \phi(C_2 u)), \ u > C,
$$
where at $ u \to \infty \ \Rightarrow \ \phi(u) \to 0. $ We denote
$$
{\bf Q}(u) = {\bf P} \left( ||\hat{f} - f||^2 / \tau^*(n) > u \right).
$$
{\bf Theorem $ \tau $.}  {\it At $ u \le C/B(n) $ the following inequality 
holds:}
$$
Q(u) \le 2 \exp( - \phi(C \sqrt{u})).
$$
{\bf Proof.} We have by the full probability formula we (we will understood 
$ {\bf P}(A \big/ B) $ 
as  the conditional probabilities, if, of course, $ A $ and $ B $ are events):
$$
{\bf Q}(u) \le {\bf P} \left( \frac{ ||\hat{f} - f||^2}{\tau^*(n)} > u \big/
\frac{\tau^*(n)}{B(n)} > \frac{1}{v}  \right) \cdot
{\bf P} \left( \frac{\tau^*(n)}{B(n)} > \frac{1}{v} \right) +
$$

$$
+ {\bf P} \left( \frac{ ||\hat{f} - f||^2}{\tau^*(n)}> u \big/ \frac{\tau^*(n)}{B(n)} \le 
\frac{1}{v} \right) \cdot {\bf P} \left( \frac{\tau^*(n)}{B(n)} \le \frac{1}{v} \right)
\stackrel{def}{=} Q_1 + Q_2;
$$

$$
Q_1 \le {\bf P} \left( || \hat{f} - f||^2/B(n) > u/v \right) 
\le \exp( - \phi(C_1 u/v));
$$

$$
Q_2 \le {\bf P} \left( \tau^*(n)/B(n) \le 1/v  \right) \le 
\exp \left( - \phi(C_2 v) \right).
$$
 Summing up and put $  v = C_3 \sqrt{u} $, we obtain  the assertion of 
the theorem. \par
 The increase in the probability $ {\bf Q} $ compared to $ {\bf P}_f $ is 
apparently explained by 
the ability of the denominator, i.e. $ \tau^*(n) $ to take values nearly to 
zero. \par
 Note in conclusion that the estimates proposed by us have successfully 
passed 
experimental tests on problems {\bf R, D } by simulate of  modeled with the 
use of pseudo-random numbers as well as on real data (of financial data) 
for which our estimations of the regression and density were compared with 
classical estimates obtained by the kernel and  wavelets
estimations method. The precision of estimations proposed here  is better. \par
 The advantage of our estimations in comparison to the trigonometrical 
estimations [Ostrovsky, Sirota, 2004] is especially in the case when 
the estimating function $ f(\cdot) $ is not periodical. \par

\vspace{2mm}

\newpage

{\bf References} \\

\vspace{3mm}

 AAD W. VAN DER VAART and MAARK J. van Der Laan. Smooth Estimation of a 
monotonic density. {\it Statistics.} 2003, V, 37, $ N^o $ 3, p. 189 - 203.\\ 

 Allal J., and Kaaouachi. Adaptive R - estimation in a Linear Regression Model with Arma 
 Errors. {\it Statistics.} 2003,  V. 37, $N^o $ 4, July - August, pp. 271 - 286.\\

 Bateman H., Erdelyi A. Higer Transcendental Functions. MC Graw-Hill Book 
Company, V.2,  New York, Toronto,…, 1953.\\

 Bobrov P.B., Ostrovsky E.I.  Confidence intervales by adaptive estimations.
{\it  Zapiski Nauchn. seminarov POMI.}  St. - Petersburg, 1997, v. 37 b.2, 28 - 45.\\ 

 Candes E.J.,  Ridgelets: Estimating with ridge Functions.  {\it Annals of 
 Statistics,} 2003, v. 31 $ N^o $ 31, 1561 - 1569.\\

 Corrine Berzin, Jose' R. Leon and Joaquim Ortega. Convergence of non - linear 
functionals of Smoothed Empirical Processes and Kernel Density Estimates. 
{\it Statistics, } 2003, V. 37 $ N^o $ 4, pp. 217 - 242. \\

 Dette H. and Melas V. Ch. B.  Optimal Design for Estimating individual 
 coefficients in Fourier Regression Model.{\it Annales of Statistics,} 2003, v. 31
 $ N^o $ 5, 1669 - 1692.\\

  DeVore R.A., Lorentz G.G. Constructive Approximation. Springer-Verlag, 
 1993.\\

 Donoho D., Jonstone I., Keryacharian G., Picard D.  Density estimation 
 by wavelet thresholding.{\it  Technical report} $ N^0 426,$  1993, Dept. of Stat., 
 Stanford University.\\ 

 Donoho D., Jonstone I.  Adapting to unknown smoothness via wavelet 
 shrinkage.  {\it Technical report} $ N^0 425,$ 1993,  Dept. of Stat.,
 Stanford University.\\

 Donoho D.  Wedgelets: nearly minimax estimation of edges.{\it  Annales of 
 of Statist.,} 1999, v. 27 b. 3 pp. 859 - 897.\\

 Donoho D.  Unconditional bases are optimal bases for data compression 
 and for statistical estimation. {\it Applied Comput. Harmon. Anal., } 1996,
 v. 3 pp. 100 - 115.\\

 Donoho D.  Unconditional bases and bit - level compression.{\it Appl. 
 Comput. Anal.,}  1999, v. 3 pp. 388 - 392.\\

 Efroimovich S.  Nonparametric estimation of the density of a unknown 
 smoothness. {\it Theory Probab. Appl.,} 1985, v. 30 b. 3,  557 - 568.\\

 Fiegel T., Hitczenko P., Jonson W.B., Shechtman G., Zinn J.  Extremal .
 properties of Rademacher functions with applications to the Khinchine and Rosental 
 inequalities. {\it Transactions of the American Math. Soc.,} 1997, v. 349 $ N^o $ 3, 
 997 - 1024.\\ 

 Golubev G., Nussbaum M.  Adaptive spline Estimations in the nonparametric 
 regression Model. {\it Theory Probab. Appl.,} 1992, v. 37 $ N^o $ 4, 521 - 529.\\

 Golybev G. {\it Nonparametric estimation of smooth spectral densities of 
 Gaussian stationary sequences.} Theory Probab. Appl., 1994, v. 38 b. 2, 
 28 - 45.\\

  Hall P., Heyde C.C. {\it Martingale Limit Theory and Applications.} Academic 
 Press,  New York, (1980) \\

 Ibragimov I.A., Khasminsky R.Z. {\it On the quality boundaries of nonparametric 
 estimation of regression.} Theory Probab. Appl., 1982, v. 21 b. 1, 81 - 94.\\

 Jonson W.B., Schechtman G., Zinn J. Best Constants in the moment Inequalities for linear
 combinations of independent and exchangeable random Variables. {\it  Annales Probab.,}
 1985, v. 13, 234 - 253.\\

 Kallaev S.O., {\it On de la Vallee Poissin sums of Fourier-Gegenbauer series}.
 Math. Zametki, {\bf 7}, $ N^o $ 1, 19-31 (1970) \\ 

 Korolyuk V.S., Borovskich Yu. P. {\it Theory of $ U - $ Statistics,} 1993, Springer,
 Berlin - Heidelberg - New York - Tokyo. \\

 Kozachenko Yu.V., Ostrovsky E.I.  Banach Spaces of random Variables of subgaussian type.A
 {\it Theory Veroyatn. Mathem. Statist.,} 1983, v. 32, 52 - 53.\\

 Lee Geunghee.  Choose of smoothing Parameters in Wavelet Series Estimators.
 Journal of Nonparametric Statistics, 2003, v. 15 (4 - 5), p. 421 - 435.\\

 Lepsky O. {\it On adaptive estimation problem in the Gaussian white Noise.} 
 Theory Probab. Appl., 1990, v. 35 b. 3, 454 - 461.\\

 Lucet Y. Faster than the Fast Legendre Transform, the Linear-time Legendre 
Transform. Numerical Algorithms, 2004, V. 16 num. 2, 171-185.\\

  Nikolsky S. {\it Inequalities for integer Functions of finite Power add Their 
 Applications in the Theory of differentiable Functions of many Variables.} (in Russian).
 In: Trudy Mathemat. Inst. im. V.V.Steklova AN SSSR, 1951, v. 51, 244 - 278.\\

 Nussbaum M. {\it Spline smoothing in regression Models and Asymptotic 
 Efficiency in } $ L_2. $  Annales of Statist., 1985, v. 13 b. 3, 984 - 997.\\

 Ostrovsky E.I. Adaptive estimation in three classical 
 problem  on nonparametric statistics. {\it Aktualnye problemy 
 sovremennoy mathematiki.}  Novosibirsk, NII MI OO, 1997, v.3 pp. 142 - 146.\\

 Ostrovsky E.I. The adaptive estimation in multidimensional statistics.A
 In: {\it Proseeding of the 5th international conference on simulation 
 of devices and technologies (ISDT).}  Obninsk, 1996, pp. 115 - 118.\\

 Ostrovsky E.I. {\it Exponential estimates for random fields and their 
 applications } (in Russian). Obninsk, OIATE, 1999, 350 pp.\\

 Ostrovsky E., Sirota L. {\it Universal adaptive estimations and confidence 
 intervals in the non-parametrical statistics.} Electronic Publications,
 arXiv.mathPR/0406535 v1 25 Jun 2004.\\

 Ostrovsky E.  {\it Bide-side exponential and moment inequalities for tails 
 of distributions of polynomial martingales.} Electronic Publications, 
 arXiv:math.PR/0406532 v1 25 Jun 2004 \\ 

 Pizier G. Condition  d'entropie  assurant la continuite de certains 
 processes et applications a  l'analyse  harmonique. In: {\it Sem. 
  d' anal. funct.,} 1979 - 1980, v. 23 - 24, 1 - 43. \\ 

 Plicusas A. Some Properties of Multiply Integral Ito.{\it Liet. Mathem. 
 Rink.,} 1981, v. 21 b. 2, 163 - 173.\\

 Polyak B., Tsybakov A. $ C_p - $ criterion in projective Estimation of 
 Regression.{\it  Theory Probab. Appl.,} 1990, v. 35 b. 2, 293 - 306.\\

 Polyak B., Tsybakov A. A family of asymptotically optimal Methods for selecting the Order of projectiv Estimation of Regression.{\it Theory Probab. Appl.,} 1992, v. 37, b. 3,  471 - 485.\\

  Potts D., Steidl G.,  Tasche M. Fast algorithms for discrete polynomial transforms. {\it Mathematics of Computation, } 1998, v. 67, n. 224, 1577-1590.\\

  Ronzin A.  Asymptotic formulae for moments of $ U \ - $ statistics with   degenerate  kernel.{\it Theory Probab. Appl.,} 1982, v. 27 b. 2, 163 - 173.\\

  Rozental H. On the subspaces of $ L_p (p > 2) $  spanned by sequence of
  independent variables. {\it Probab. Theory Appl.,} 1982, v. 27 b.1, 47 - 55.\\
 
  Saulis L., Statuliavichius V. {\it Limit Theorems for Great Deviations.} 
  Vilnius, Mokslas, (1989) (in Russian) \\

  Shiryaev A.N. {\it Probability.} Kluvner Verlag, 1986.\\ 

  Szeg\"o G. {\it Orthogonal Polynomials.} Amer. Math. Soc. Colloq. Publ.,
  Vol. 23, New York, (1959).\\
   
  Tchentsov N.N. {\it Statistical decision rules and optimal inference.}
  Moscow, Nauka, 1972. \\

  Timan A. {\it Theory of Approximation of Functions of Real Variables } (in Russian).
  Moscow, GIFML, 1960.\\
   
  Tony Cai T.  Adaptive wavelet Estimation: a Block Thresholding and Oracle 
  Inequality Approach.{\it  Annales of Math. Statist.,} 1999, v. 27 b. 3, 898 - 924.\\

\end{document}